\documentclass[12pt]{amsart}
\usepackage{amsfonts, amssymb, amsmath, amsthm, latexsym, array}
\usepackage{fullpage}
\usepackage{verbatim,euscript}

\input {cyracc.def}
 \font\tencyr=wncyr10 
\def\rus{\tencyr\cyracc}

\newtheorem{thm}{Theorem}[section]
\newtheorem{lm}[thm]{Lemma}
\newtheorem{cl}[thm]{Corollary}
\newtheorem{prop}[thm]{Proposition}

\theoremstyle{remark}

\theoremstyle{definition}
\newtheorem{ex}[thm]{Example}
\newtheorem{df}{Definition}

\newcommand {\ce}{{\mathfrak c}}

\newcommand {\g}{{\mathfrak g}}

\newcommand {\es}{{\mathfrak s}}

\newcommand{\um}{{\cdot}}
\newcommand{\umt}{{\times}}

\newcommand {\esi}{\varepsilon}
\newcommand{\lb}{\lambda}
\newcommand{\ap}{\alpha}
\newcommand{\vp}{\varphi}

\renewcommand{\le}{\leqslant}
\renewcommand{\ge}{\geqslant}

\def\Hom{\mathrm{Hom}}
\def\ad{\mathrm{ad}\,}

\def\Ad{\mathrm{Ad}}

\def\SO{{\rm SO}}
\newcommand{\Spin}{{\rm{Spin}}}
\newcommand {\pn}{{\it pn}-}
\newcommand{\odin}{{\mathrm{1\hspace{1pt}\!\! I}}}

\newcommand{\gt}{\mathfrak}
\newcommand{\eus}{\EuScript}
\newcommand{\BB}{\mathbb}
\newcommand{\id}{{\rm id}}
\newcommand{\GL}{{\rm GL}}
\newcommand{\cgo}{\gt C(\g_1)}

\newcommand{\Lie}{{\mathrm{Lie\,}}}

\newcommand {\rk}{\mathrm{rk\,}}

\newcommand {\Ker}{{\mathrm{Ker\,}}}
\newcommand {\fN}{\mathfrak{N}}

\newfam\eusfam%
\font\Bbbfont=msbm10 scaled 1200%
\def\bbk{\hbox {\Bbbfont\char'174}}

\begin{document}
\setlength{\parskip}{3pt plus 5pt minus 0pt}
\hfill {\scriptsize January 22, 2006}
\vskip1ex

\title[Symmetric pairs and commuting varieties]
{Symmetric pairs and associated commuting varieties}
\author[D.\,Panyushev]{Dmitri Panyushev}
\address[D.\,Panyushev]{Independent University of Moscow,
Bol'shoi Vlasevskii per. 11, 119002 Moscow, \ Russia}
\email{panyush@mccme.ru}
\author[O.\,Yakimova]{Oksana Yakimova}
\address[O.\,Yakimova]{MPI f\"ur Mathematik,
Vivatsgasse 7, \ 53111 Bonn, Germany}
\email{yakimova@mpim-bonn.mpg.de}
\thanks{This research was supported in part by
CRDF Grant RM1-2543-MO-03 and RFBI Grant 05-01-00988}
\maketitle

\section*{Introduction}

The ground field $\bbk$ is algebraically closed and of characteristic zero.
Let $\g$ be a reductive algebraic Lie algebra over $\bbk$
and $\sigma$ an involutory automorphism of $\g$.
Then $\g=\g_0\oplus\g_1$ is the direct sum of $\sigma$-eigenspaces.
Here $\g_0$ is a reductive subalgebra and
$\g_1$ is a $\g_0$-module. We also say that $(\g,\g_0)$ is
a {\it symmetric pair}.
Let $G$ be the adjoint group of $\g$
and $G_0$ the connected subgroup of $G$ with $\Lie G_0=\g_0$.
The {\it commuting variety} associated with the 
involution $\sigma$ or symmetric pair 
$(\mathfrak{g},\g_0)$  is:
\[
{\gt C}(\g_1) = \{ (x,y) \in \g_1 \times \g_1 \ \mid \ [x,y]=0\}.
\]
The "usual" commuting variety is obtained as a special case of this 
construction. Indeed, let $(\g\oplus\g,\g)$ 
be the symmetric pair such that $\g$ 
is embedded diagonally into $\tilde{\g}=\g\oplus\g$. Then 
$\tilde{\g}_1\simeq\g$ as $\g$-module and
$\gt C(\tilde{\g}_1)$ is isomorphic to 
the commuting variety ${\gt C}(\g)=\{(x,y)\in\g\umt
\g\mid[x,y]=0\}$. The latter was
considered by  Richardson in \cite{r}. He proved that
${\gt C}(\g)=\overline{G(\gt t{\times}\gt t)}$, where $\gt t\subset\g$
is a Cartan subalgebra, and therefore ${\gt C}(\g)$ is irreducible.
It is not yet known whether ${\gt C}(\g)$ is normal or whether the 
ideal of ${\gt C}(\g)$ in $\bbk[\g\times\g]$
is generated by quadrics.

Let $\mathfrak{c}\subset\g_1$ be a maximal abelian subspace
consisting of semisimple elements. Any such subspace is
called a {\it Cartan subspace} of $\g_1$. All Cartan
subspaces are $G_0$-conjugate, see \cite{kr}. 
The {\it rank} of the symmetric pair
$(\gt{g}, \g_0)$ is $\dim\gt c$, also denoted $\rk(\g,\g_0)$.
It is known that $\dim{\gt C}(\g_1)=\dim\g_1 + \dim\gt{c}$ 
and ${\gt C}_{0}:=
\overline{G_0(\gt{c}\umt  \gt{c})}$ is an
irreducible component of ${\gt C}(\g_1)$ of maximal
dimension, see \cite[Sect.\,3]{Dima1}. It follows that $\gt C(\g_1)$ is 
irreducible if and only if $\gt C(\g_1)=\gt C_0$.

The irreducibility problem for $\gt C(\g_1)$ was first considered 
by Panyushev \cite{Dima1}. 
He noticed that $\gt C(\g_1)$ can be
reducible \cite[3.5]{Dima1}. On the other hand, 
in some particular cases $\gt C(\g_1)$ is irreducible and enjoys
a number of good properties.
If $(\g,\g_0)$ is a symmetric pair of
maximal rank (i.e., $\rk(\g,\g_0)=\rk\g$), 
then $\gt C(\g_1)$ is an irreducible normal
complete intersection and the ideal of $\gt C(\g_1)$
in $\bbk[\g_1{\times}\g_1]$
is generated by quadrics, see \cite[3.2]{Dima1}.
We refer to \cite{Dima2},\,\cite{sy},\,\cite{sy2} for some recent results.

The problem of irreducibility
is closely related to some properties of rings of differential operators. 
Let ${\mathcal D}(\g_1)$  denote the algebra of 
differential operators on $\g_1$ with coefficients in $\bbk[\g_1]$
and ${\mathcal K}\lhd {\mathcal D}(\g_1)$ the left ideal annihilating 
all $G_0$-invariant polynomials.
The $G_0$-action on $\g_1$ allows us to interpret the
elements of $\g_0$ as left-invariant vector fields on $\g_1$.
This gives rise to a homomorphism of the 
universal enveloping algebra ${\bf U}(\g_0)$ 
to ${\mathcal D}(\g_1)$. If $\g_0$ contains no simple ideals of
$\g$, then this homomorphism is injective. We may
therefore assume that ${\bf U}(\g_0)\subset{\mathcal D}(\g_1)$. 
It is easily seen that ${\bf U}(\g_0)\subset{\mathcal K}$. Furthermore,
it can be shown that if  
${\gt C}(\g_1)$ 
is irreducible and its ideal is generated by quadrics, then
${\mathcal K}={\mathcal D}(\g_1){\bf U}(\g_0)$. 
(See \cite{LSt} for related results.)

\vskip.8ex
Suppose $\g$ is simple. Then the known results are
\begin{itemize}
\item[$\bullet$] if 
$\rk(\g,\g_0)=\rk\g$ 
(the {\it maximal rank} case), then
${\gt C}(\g_1)$ is irreducible \cite{Dima1};
\item[$\bullet$]  if $\rk(\g,\g_0)=1$, then ${\gt C}(\g_1)$ 
is reducible unless
$(\g,\g_0)=(\gt{so}_{m+1}, \gt{so}_{m})$  \cite{sy}; 
furthermore, the number of irreducible components of $\cgo$ equals
the number of non-zero nilpotent $G_0$-orbits in $\g_1$ \cite{Dima2};
\item[$\bullet$] for
$(\gt{sl}_{2n}, \gt{sp}_{2n})$ and $({\EuScript E}_6,{\EuScript F}_4)$, 
the commuting variety $\cgo$ is irreducible \cite{Dima2};
\item[$\bullet$] if $(\g,\g_0)=(\gt{so}_{m+2},
\gt{so}_2{\oplus}\gt{so}_m)$, then ${\gt C}(\g_1)$
is irreducible \cite{sy2}.
\end{itemize}
The goal of this article is to present some further results on the
irreducibility problem.
Let $\g_x$ denote the centraliser of $x\in\g$. If $x\in\g_1$, 
then $\sigma$ induces the symmetric decomposition 
$\g_x=\g_{0,x}\oplus\g_{1,x}$, where $\g_{0,x}:=\g_0\cap\g_x$ is the
centraliser of $x$ in $\g_0$.
If $h\in\g_1$ is semisimple, then $(\g_h,\g_{0,h})$ is a symmetric
pair, which is called a {\it sub-symmetric pair\/} associated with 
$(\g,\g_0)$.
Let $\gt C(\g_{1,h})$ denote the commuting variety associated with 
$(\g_h,\g_{0,h})$. Usually, the irreducibility 
of $\gt C(\g_1)$ is proved in the following way, which goes back 
essentially to Richardson \cite{r}.  
First one has to show that $\gt C(\g_{1,h})$ is irreducible for any 
non-zero semisimple $h\in\g_1$; 
next, one has to verify that 
$(\{e\}{\times}\g_{1,e})\subset\gt C_0$ for any
nilpotent element $e\in\g_1$. Actually, 
the verification in case of nilpotent elements 
readily reduces to $\sigma$-distinguished elements (see Section~\ref{1}
for precise definitions).
To a great extent, the structure of $\gt C(\g_1)$ depends on
properties of $\sigma$-distinguished elements.

In Section~\ref{1}, we collect several useful facts concerning semisimple 
and nilpotent elements in $\g_1$. To implement the above program, we should
be able to deal with the respective sub-symmetric pairs, and 
Proposition~\ref{sat} provides a description of all possible sub-symmetric pairs 
$(\g_h,\g_{0,h})$ in terms of the Satake diagram of $(\g,\g_0)$. 

In Section~\ref{f}, we prove the irreducibility of the commuting variety 
for $(\gt{so}_{n+m},\gt{so}_n\oplus\gt{so}_m)$.
This extends the result of \cite{sy2}, which refers to the case $n=2$.
The scheme of the proof is similar to that of \cite{sy2}. But as it often 
happens, the argument in a general situation is shorter and simpler than
in a particular case. Roughly speaking, the reason of success is that,
for this symmetric pair, all $\sigma$-distinguished elements are even.

In Section~\ref{e6}, $\gt C(\g_1)$ is shown to be
irreducible for the symmetric pairs $(\gt{sl}_{2n},
\gt{sl}_n{\oplus}\gt{sl}_n{\oplus}{\gt t}_1)$ and $({\EuScript E}_6,
\gt{sl}_6{\oplus}\gt{sl}_2)$. 
Here we use the fact
that in these two cases $\dim\g_1>\dim\g_0$ and the (closures of) 
subsets $G_0(\{e\}{\times}\g_{1,e})$ cannot form an irreducible  
component of $\gt C(\g_1)$. In the classical case, 
$\gt C(\g_{1,h})$ is shown to be irreducible by induction, 
and in the ${\EuScript E}_6$-case we rely on the explicit description of 
all sub-symmetric pairs $(\g_h,\g_{0,h})$.

In \cite{Dima2}, it was conjectured that ${\gt C}(\g_1)$
is irreducible if $\rk(\g,\g_0)>1$. This is, however, false. 
In Section~\ref{nex},
we prove that $\gt C(\g_1)$ is reducible for the symmetric
pairs $(\gt{sl}_{n+m},\gt{sl}_n{\oplus}\gt{sl}_m{\oplus} {\gt t}_1)$ with
$n\ne m$, $(\gt{so}_{2n},\gt{gl}_n)$ with odd $n$, and
$({\eus E}_6, \gt{so}_{10}\oplus {\gt t}_1)$.
Their ranks are equal to $\min(n,m)$, $[n/2]$, and $2$, respectively.
Furthermore, we give am explicit lower bound on the number of the
irreducible components of $\gt C(\g_1)$ in the $\gt{sl}_{n+m}$-case,
which shows that this number can be arbitrarily large, see 
Proposition~\ref{irrcomp}. In Section~\ref{pnp}, we prove the reducibility
of $\cgo$ for $({\eus E}_7,\gt{so}_{12}{\oplus}\gt{sl}_2)$ 
and $({\eus E}_8,{\eus E}_7{\oplus}\gt{sl}_2)$, using some ideas form
the theory of principal nilpotent pairs.

Taking into account the previously known results and the results of this paper,
one sees that there remain only three cases (two classical series and one
symmetric pair for $\eus E_7$), where the the answer is not known.
In particular, for series $\eus A$ and $\eus B$,
the irreducibility problem is completely  solved. 
In Section~\ref{unk}, we discuss possible approaches
to the remaining cases. 

\noindent
{\bf Acknowledgements.} {\small This paper was written during 
our stay at the Max-Planck-Institut f\"ur Mathematik (Bonn). 
We are grateful to this institution for the warm hospitality and support.}

\section{Nilpotent elements, sub-symmetric pairs, and commuting varieties}
\label{1}

In this section, we deal with semisimple and 
nilpotent elements in $\g_1$ and their relations to 
the commuting variety. 

Let $n$ be a non-negative integer. The set
\[
\g_1^{(n)}=\{\xi\in \g_1\mid \dim(G_0\xi)=n \}
\]
is locally closed. The irreducible components of the sets $\g_1^{(n)}$
are called the {\it $G_0$-sheets} of $\g_1$. 

\begin{lm}\label{sheet}
Let $\eus S$ be a $G_0$-sheet of $\g_1$ containing semisimple
elements. Suppose that, for each semisimple $h\in \eus S$, we have
$(\{h\}\umt  \g_{1,h})\subset\gt C_0$. Then $(\{x\}\umt
\g_{1,x})\subset\gt C_0$ for each $x\in \eus S$.
\end{lm}
\begin{proof}
Let $x\in \eus S$. Since $\eus S$ contains semisimple
elements, they form a dense open subset. Therefore,
we can find a morphism
$\gamma:\bbk\to \eus S$ such that
$\gamma(0)=x$ and $\gamma(t)$ is semisimple
for each $t\ne 0$. Then
$\g_{1,x}=\lim_{t\to 0}\g_{1,\gamma(t)}$, where the
limit is taken in an appropriate Grassmannian.
For each $y\in\g_{1,x}$, we can define
elements $y(t)\in\g_{1,\gamma(t)}$ such that
$y=\lim_{t\to 0}y(t)$.
Since $(x,y)=\lim_{t\to 0}(\gamma(t),y(t))$ and
$(\gamma(t),y(t))\in\gt C_0$ for each $t\ne 0$,
we conclude that $(x,y)\in\gt C_0$.
\end{proof}

\begin{lm}\cite[Prop.\,5]{kr}\label{z2} 
For any $x\in\g_1$, we have $\dim\g_{1,x}-\dim\g_{0,x}=\dim\gt
g_1-\dim\g_0$.
\end{lm}

Let $\fN(\g_1)$ denote  
the set of all nilpotent elements in $\g_1$.
For any $e\in\fN(\g_1)$, there is an $\gt{sl}_2$-triple
$\{e,f,h\}$ such that $f\in\g_1$ and $h\in\g_0$ \cite{kr}.
Such a triple is said to be {\it normal}.
Recall that $e$ is called {\it even} if the eigenvalues of
$\ad h$ on $\g$ are even. 

\begin{df}
An element $e\in\fN(\g_1)$
is said to be {\it $\sigma$-distinguished} (in other
notation, $\gt p$- or $({-}1)$-distinguished) if $\g_{1,e}$ contains 
no semisimple elements of $[\g,\g]$.
\end{df}

The following lemma appears in \cite{sy2}, but the proof
given here is shorter.

\begin{lm}\label{even}
Suppose $e\in\fN(\g_1)$ is even. Then
$e$ belongs to a $G_0$-sheet containing semisimple elements.
\end{lm}
\begin{proof}
Let $(e,f,h)$ be a normal $\gt{sl}_{2}$-triple.
Since $e$ is even, we have $\dim\gt
g_h=\dim\g_e$. Set $e(t):=e-t^2f$  for $t\in\bbk$. If
$t\ne 0$, then $e(t)$ is semisimple and conjugate to $th$.
Therefore $\dim\g_{e(t)}=\dim\g_h=\dim\g_e$ and 
$\dim\g_{0,e(t)}=\dim\g_{0,e}$  by
Lemma~\ref{z2}.
Clearly $e(0)=e=\lim_{t\to 0} e(t)$ and the $G_0$-sheet
containing $e$ contains also semisimple elements $e(t)$.
\end{proof}

Let $h\in\g_1$ be a semisimple element. Recall that $\gt
C(\g_{1,h})$ denotes for the commuting variety associated with the
sub-symmetric pair $(\g_h,\g_{0,h})$.

\begin{thm}\label{main}
Suppose that each $\sigma$-distinguished nilpotent element
in $\g_1$ is even and ${\gt C}(\g_{1,h})$
is irreducible for each semisimple non-zero $h\in g_1$.
Then ${\gt C}(\g_1)$ is irreducible.
\end{thm}
\begin{proof}
Recall that ${\gt C}_{0} =\overline{G_0(\gt c\umt  \gt c)}$
is an irreducible component of $\cgo$. Following the
original proof of Richardson \cite{r} (see also
\cite[Sect.~2]{Dima2}), we show that $\gt C(\g_1)=\gt
C_0$. Let $(x,y) \in {\gt C}(\g_1)$.

(1) Suppose
there is a  semisimple element $h\in\g_1$ such that $[h,x]=[h,y]=0$.
This assumption is automatically satisfied if either $x$ or $y$
is semisimple. Moreover, if $x$ (or $y$) is not nilpotent and
$x=x_s+x_n$ is the Jordan decomposition, then
$x_s\in\g_1$ and $[x_s,x]=0$, $[x_s,y]=0$.

Consider the sub-symmetric pair $(\g_h, \g_{0,h})$.
Replacing $\gt c$ with a conjugate Cartan subspace, we may
assume that $h\in\gt c$. Then $\gt c$ is a Cartan subspace
of $\g_{1,h}$. By the assumption, ${\gt C}(\g_{1,h})$ 
is irreducible. 
Therefore, $(x,y)\in\overline{G_{0,h}(\gt c\umt
\gt c)}$ and hence $(x,y)\in\gt C_0$.

(2) Now we may assume that both $x$ and $y$ are nilpotent.
Suppose first that there is a semisimple element
$h\in\g_1$ such that $[x,h]=0$.
Then $(x,(1-t)y+th)\in\gt C(\g_1)$ for each
$t\in\bbk$ and $(1-t)y+th$ is nilpotent only for
a finite number of $t$'s. Therefore, by part~(1),
one has $(x,(1-t)y+th)\in\gt C_0$ for almost all $t$.
Since $y=\lim_{t\to 0}((1-t)y+th)$, we get $(x,y)\in\gt C_0$.

(3) It remains to handle the case in which both $x$ and $y$ are
$\sigma$-distinguished nilpotent elements. Then both $x$ and $y$ 
are even and by Lemma~\ref{even}, $x$ belongs to a
$G_0$-sheet containing semisimple elements.
According to part~(1), the assumptions of
Lemma~\ref{sheet} are satisfied and it follows that
$(x,y)\in\gt C_0$.
\end{proof}
In order to effectively apply Theorem~\ref{main}, we need a 
description of possible sub-symmetric pairs.
To this end, we recall the notion of the 
{\it Satake diagram} of a symmetric pair $(\g,\g_0)$, where $\g$ is
semisimple. Roughly speaking, it 
is the Dynkin diagram of $\g$, where each node is either black or white
and some pairs of white nodes are connected by a (new) arrow.
More precisely, let $\gt t$ be a $\sigma$-stable Cartan subalgebra of
$\g$ containing $\gt c$ and let $\Delta$ be the root system of $(\g,\gt t)$.
Since $\gt t$ is $\sigma$-stable, $\sigma$ acts on $\Delta$.
It is possible to choose the set of positive roots, $\Delta^+$, such that
if $\beta\in\Delta^+$ and $\beta(\gt c)\ne 0$, then 
$\sigma(\beta)\in -\Delta^+$. Let $\Pi$ be the corresponding set of 
simple roots and $\ap\in\Pi$.
Then the node corresponding to $\ap$ is black if and 
only if $\alpha(\gt c)=0$.
Two different white nodes corresponding to 
$\alpha_i,\alpha_j\in\Pi$ are connected by an arrow if
and only if $\alpha_i\vert_{\gt c}=\alpha_j\vert_{\gt c}$.
The latter is equivalent to that $\sigma(\ap_i)=-\ap_j$.
Looking at the Satake diagram, one immediately reads off
many properties of a symmetric pair under consideration. For instance,
\[
\rk(\g,\g_0)=\text{(number of white nodes)}-\text{(number of arrows)},
\]
and the centraliser of a generic element of $\ce$ is the Levi subalgebra
of $\g$ whose semisimple part is given by the subdiagram consisting
of all black nodes.
\\
Let $S$ be a Satake diagram. Then $S'$ is a {\it subdiagram of\/} $S$
if $S'$ is obtained from $S$ by repeating the following 
procedure. Namely, one can
remove one white node, if it is not connected by an arrow; or one 
can remove a pair of nodes connected by an arrow.
\\
Let $h\in\g_1$ be semisimple. Set $\gt p=[\g_h,\g_h]$ and
$\gt p_i=\g_i\cap\gt p$. Then $(\gt p,\gt p_0)$ is said to be a
{\it reduced\/} sub-symmetric pair of $(\g,\g_0)$.
Clearly, $\gt C(\g_{1,h})$ is irreducible if and only if
$\gt C(\gt p_1)$ is irreducible. So, it is enough to describe all
reduced sub-symmetric pairs.

\begin{prop}     \label{sat} 
A symmetric pair $(\gt p, \gt p_0)$ is a reduced sub-symmetric pair of
$(\g,\g_0)$ if and only if $S_{\gt p}$, the Satake diagram  
of $(\gt p, \gt p_0)$, 
is a sub-diagram of $S_{\g}$, the Satake diagram of $(\g,\g_0)$.
\end{prop}
\begin{proof} Without loss of generality, 
we may assume that $h\in\gt c$. Then $\ap(h)=0$ if the node 
corresponding to $\ap\in\Pi$ is black. Since we are only interested
in possible types of centralisers, we may assume that $h$ lies in the $\BB Q$-form
$\ce_{\BB Q}$. Here $\ce_{\BB Q}$ is the $\BB Q$-span of all restricted roots
$\beta\vert_\ce$, $\beta\in\Delta$, and $\ce$ is identified with $\ce^*$ using the 
restriction of the Killing form.
Moreover, using the action of the little Weyl group on $\ce_{\BB Q}$, 
we may assume that $\ap(h) > 0$ if the node corresponding to
$\ap\in\Pi$ is white. Then the Dynkin diagram $D_{\gt p}$
of $\gt p$ is a subdiagram of 
$D_\g$ consisting of all $\alpha\in\Pi$ such that
$\alpha(h)=0$. In particular, $D_{\gt p}$ contains all black nodes
of $S_\g$.

Set $\gt c_{\gt p}:=\gt c\cap\gt p$. It is a Cartan subspace of $\gt p_1$
and $\gt c=\gt c_{\gt p}\oplus\gt z_1$,
where $\gt z_1$ lies in the centre of $\g_h$.
If $\alpha_i$ is a simple root of $\gt p$, then
$\alpha_i(\gt z_1)=0$ and $\alpha_i(\gt c_{\gt p})=0$
if and only if $\alpha_i(\gt c)=0$. Hence
the nodes do not change their colours.
Suppose $\alpha_i$ and $\alpha_j$ are connected by an arrow
in $S_{\g}$. Then $(\ap_i-\ap_j)(\gt c)=0$. In particular,
$\alpha_i(h)=0$ if and only if $\alpha_j(h)=0$.
If $\alpha_i$ and $\alpha_j$ belong to $D_{\gt p}$,
then they are still connected by an arrow, since 
$\gt c_{\gt p}\subset\gt c$.
Conversely, if $\alpha_i,\alpha_j$ are simple roots of 
$\gt p$ and they are connected by an arrow of $S_{\gt p}$, 
then $(\alpha_i-\alpha_j)(\gt c_{\gt p})=0$ and 
$(\alpha_i-\alpha_j)(\gt z_1)=0$. 
Therefore $(\alpha_i-\alpha_j)(\gt c)=0$, and 
$\alpha_i$, $\alpha_j$ are connected by an arrow in
$S_{\g}$. Thus, the Satake diagram of any reduced sub-symmetric pair is
a subdiagram of $S_\g$.

Conversely, let $S'$ be a subdiagram of $S_\g$ in the sense of
the above definition. Consider the subspace  
$\gt c'\subset\gt c$, where all simple roots corresponding to the nodes
from $S_\g\setminus S'$ vanish.
Then, for a generic $h\in\gt c'$, we obtain the reduced sub-symmetric pair
whose Satake diagram is $S'$.
\end{proof}

\section{${\gt C}(\g_1)$ is irreducible for
$(\g,\g_0)=(\gt{so}_{n+m},\gt{so}_n{\oplus}\gt{so}_m)$} 
\label{f}

In this section $G=\SO_{n+m}$, 
$\g=\gt{so}_{n+m}$, and $\g_0=\gt{so}_n{\oplus}\gt{so}_m$.
Let $V=\bbk^{n+m}$ be
a vector space of
the defining representation of $\g$.
Then we have a $G_0$-invariant decomposition
$V=V_a{\oplus}V_b$, where
$V_a=\bbk^n$, $V_b=\bbk^m$ and
$\g_1\cong\bbk^n{\otimes}\bbk^m$ as a $G_0$-module.
Let $(\,\,,\,)$  denote  the non-degenerate symmetric
$G$-invariant bilinear form on $V$.
In order to apply Theorem~\ref{main},
we need a description of semisimple and nilpotent
elements in $\g_1$.

\begin{lm}\label{ss}
Let $h\in\g_1$ be a semisimple element.
Then the symmetric pair $(\g_h, \g_{0,h})$ is a direct sum
$(\bigoplus_{i=1}^{r}(\gt{gl}_{k_i},\gt{so}_{k_i}))\oplus
(\gt{so}_{n+m-2k},\gt{so}_{n-k}{\oplus}\gt{so}_{m-k})$,
where $k=\sum_i k_i$.
\end{lm}
\begin{proof}
Let $v_\lambda\in V$ be an eigenvector of $h$ such that
$h\cdot v_\lambda=\lambda v_\lambda$ and $\lambda\ne 0$.
Since $h$ preserves the symmetric form
$(\,\,,\,)$, i.e., $(h{\cdot}v_\lb,v_\lb)+(v_\lb,h{\cdot}v_\lb)=0$,
we have $(v_\lambda, v_\lambda)=0$.
Also, if $h{\cdot}v=\lambda v$, $h{\cdot}w=\mu w$, then
$(v,w)=0$ unless $\lambda=-\mu$.
Let $\{\pm\lambda_i, 0\mid i=1,\ldots, r\}$
be the set of the eigenvalues of $h$.
Then there is an orthogonal
$h$-invariant decomposition
$$
V=(V_{\lambda_1}{\oplus}V_{-\lambda_1})\oplus\ldots\oplus
(V_{\lambda_r}{\oplus}V_{-\lambda_r})\oplus V_0.
$$
Here each $V_{\lambda_i}$ is an isotropic subspace,
$(V_{\lambda_i}, V_{\lambda_j})=0$ unless
$\lambda_i+\lambda_j=0$, and $(V_0,V_{\pm\lambda_i})=0$
for each $\lambda_i$. Therefore
$\g_h\subset(\bigoplus_{i=1}^r \gt{so}(V_{\lambda_i}{\oplus}V_{-\lambda_i}))
\oplus\gt{so}(V_0)$. More precisely,
if $\dim V_{\lambda_i}=k_i$ and $k=\sum_i k_i$, then
$\g_h=(\bigoplus_i\gt{gl}_{k_i})\oplus\gt{so}_{n+m-2k}$.

Now it remains to describe $\g_{0,h}=(\g_h)^\sigma$.
We may assume that $\sigma$ is a conjugation by
a diagonal matrix $A\in{\rm O}_{n+m}$
such that $A|_{V_a}=-\id$ and $A|_{V_b}=\id$.
Since $\sigma(h)=-h$, we have
$A\cdot V_{\lambda_i}=V_{-\lambda_i}$ and $A\cdot V_0=V_0$.
Moreover,
$A$ determines a non-degenerate symmetric form $(\,\,,\,)_A$
on each $V_{\lambda_i}$ by the formula
$(v,w)_A=(v,A\cdot w)$.
Therefore, each $\gt{so}(V_{\lambda_i}{\oplus}V_{-\lambda_i})$
is $\sigma$-invariant,
$(\gt{so}(V_{\lambda_i}{\oplus}V_{-\lambda_i}))^\sigma\cong
\gt{so}_{k_i}{\oplus}\gt{so}_{k_i}$, and
$(\gt{gl}_{k_i})^\sigma=\gt{so}_{k_i}$. Finally,
the restriction $A|_{V_0}$ has signature $(n-k,m-k)$.
Thus $(\gt{so}_{n+m-2k})^\sigma=\gt{so}_{n-k}{\oplus}\gt{so}_{m-k}$.
\end{proof}
\noindent
Note that this result can also be deduced from Proposition~\ref{sat}.

Recall several standard facts concerning nilpotent elements
in $\gt{gl}(V)$. Suppose $e\in\gt{gl}(V)$ is nilpotent and
$m=\dim \Ker(e)$. By the theory of Jordan normal form,
there are vectors $w_1,\ldots,w_m\in V$ and non-negative
integers $d_1,\ldots,d_m$ such that $e^{d_i}\um w_i=0$ and
$\{ e^s\um w_i \mid 1\le i\le m, \ 0\le s <d_i\}$ is a
basis for $V$. Let $V[i]\subset V$ be the linear span of
$\{w_i,e\um w_i,\ldots, e^{d_i-1}\um w_i\}$. The spaces
$\{V[i]\}$ are called the Jordan (or {\it cyclic}) spaces of
$e$, and $V=\oplus_{i=1}^m V[i]$.

\begin{lm}     \label{ab}
Suppose $e\in\fN(\g_1)$.  Then
the cyclic vectors $\{w_i\}_{i=1}^m$ and hence the cyclic spaces
$\{V[i]\}$'s can be chosen such that the following properties are 
satisfied:
\begin{itemize}
\item[\sf (i)] \
there is an involution $i\mapsto i^*$ on the set
$\{1,\dots, m\}$ such that:  $d_i=d_{i^*}$; \ 
$i=i^*$ if and only if $\dim V[i]$ is odd; \  
$(V[i], V[j])=0$ if $i\ne j^*$;
\item[\sf (ii)] \ $\sigma(w_i)=\pm w_i$.
\end{itemize}
\end{lm}
\begin{proof}
Part (i) is a standard property of the nilpotent orbits in
$\gt{so}(V)$, see, for example, \cite[Sect.~5.1]{cm} or \cite[Sect.~1]{ja}.
Then part (ii) says that in the presence of the involution $\sigma$
cyclic vectors for $e\in\fN(\g_1)$ can be chosen to be $\sigma$-eigenvectors,
see \cite[Prop.~2]{ohta}.
\end{proof}
For each $e\in\fN(\g_1)$, we choose cyclic vectors
$\{w_i\}$ as prescribed by Lemma~\ref{ab}. Let us say that $e^s\um
w_i$ is of type $a$ if $\sigma(e^s\um w_i)=-e^s\um w_i$, i.e.,
$e^s\um w_i\in V_a$; and  $e^s\um w_i$ is of type $b$ if
$e^s\um w_i\in V_b$. Since $\sigma(e)=-e$, if $e^s\um w_i\in
V_a$, then $e^{s+1}\um w_i\in V_b$ and vice versa. Therefore
each string $( w_i,e\um w_i,\ldots,e^{d_i-1}\um w_i)$
has one of the following types: \\
\phantom{$n$} \qquad\quad
$aba\ldots ab$, \qquad  $bab\ldots ba$, \qquad $aba\ldots ba$,
\qquad $bab\ldots ab$.

\begin{prop}         \label{ev}
Every $\sigma$-distinguished element $e\in\fN(\g_1)$ 
is even.
\end{prop}
\begin{proof}
Suppose $e\in\fN(\g_1)$ is $\sigma$-distinguished
and the cyclic spaces
$\{V[i]\}$ are chosen as prescribed by Lemma~\ref{ab}.
Assume that $\dim V[i]$ is even for some $i$.
By \cite[Prop.~2]{ohta}, if
$V[i]$ is of type $aba\ldots ab$, then
$V[i^*]$ is necessarily of type $bab\ldots ba$, and vice versa;
i.e., if $\sigma(w_i)=\mp w_i$, then $\sigma(w_{i^*})=\pm w_{i^*}$.
It is easily seen that $\gt{so}(V[i]{\oplus}V[i^*])\cap\g_e$
contains a $\sigma$-stable subalgebra $\gt l$ that is isomorphic to
$\gt{sl}(\bbk w_i\oplus \bbk w_{i^*})\cong\gt{sl}_2$.
More precisely, any unimodular transformation of 
$\bbk w_i\oplus \bbk w_{i^*}$ can be uniquely extended to an element
of $\gt{so}(V[i]{\oplus}V[i^*])\cap\g_e$. 
Since the restriction of $\sigma$ to $\gt l$ is non-trivial,
$\gt l_1=\gt l\cap\g_1\subset \g_e$ contains semisimple
elements. This means that such $e$ is not $\sigma$-distinguished. Hence, 
all $V[i]$ are odd-dimensional. Since a nilpotent element of $\gt{so}_N$ 
is even if and only if all the numbers $\dim V[i]$ have the same parity, 
we are done.
\end{proof}

\begin{thm}   \label{reduction}
For $(\g,\g_0)=(\gt{so}_{n+m}, \gt{so}_n{\oplus}\gt{so}_m)$,
the commuting variety ${\gt C}(\g_1)$ is irreducible.
\end{thm}
\begin{proof}
Let $\gt c\subset\g_1$ be a Cartan subspace. We argue by
induction on $\dim\gt c$. The base of induction is the rank~$1$
case $(\gt{so}_{n+1},\gt{so}_n)$, where the irreducibility
of $\gt C(\g_1)$ is proved in \cite{Dima2}, \cite{sy}.

Let $h\ne 0$ be a semisimple element of $\g_1$.
By Lemma~\ref{ss}, $(\g_h, \g_{0,h})=
(\bigoplus_{i=1}^{r}(\gt{gl}_{k_i},\gt{so}_{k_i}))\oplus
(\gt{so}_{n+m-2k},\gt{so}_{n-k}{\oplus}\gt{so}_{m-k})$.
Each pair $(\gt{gl}_{k_i},\gt{so}_{k_i})$ is a symmetric
pair of maximal rank, hence, the corresponding commuting
variety is irreducible, see \cite[(3.5)(1)]{Dima1}. Clearly,
the commuting variety corresponding to a direct sum of
symmetric pairs is the product of the commuting
varieties corresponding to the summands. Therefore, using
the inductive hypothesis, we conclude that $\gt C(\g_{1,h})$ is 
irreducible. Taking into account Proposition~\ref{ev},
we see that the assumptions of Theorem~\ref{main}
are satisfied and $\gt C(\g_1)$ is irreducible.
\end{proof}

\section{Irreducibility of ${\gt C}(\g_1)$ for symmetric pairs with
$\dim\g_1>\dim\g_0$}\label{e6}

In this section we show that the method of proving irreducibility of
$\gt C(\g_1)$ in the maximal rank case \cite{Dima1}
actually applies in a more general setting of symmetric pairs with
$\dim\g_1>\dim\g_0$. 
Practically, this yields the irreducibility for two new pairs:
$(\gt{sl}_{2n},\gt{sl}_n{\oplus}\gt{sl}_n{\oplus} {\gt t}_1)$ and
$({\EuScript E}_6,\gt{sl}_6{\oplus}\gt{sl}_2)$.
\\[.6ex]
Recall that $(\g,\g_0)$ is of maximal rank if and only if the Satake 
diagram has no new arrows and all the nodes are white.

\begin{lm}\label{dim} Suppose all nodes of the Satake diagram
of $(\g,\g_0)$ are white and there is at least one
that is not connected by an arrow with any of the others.
Then $\dim\g_1>\dim\g_0$.
\end{lm}
\begin{proof} Choose a semisimple element $h\in\g_1$ such
that $\g_{1,h}=\gt c$ is a Cartan subspace. Then $(\g_0)_h$ 
coincides with the centraliser of $\gt c$ in $\gt
g_0$. Therefore the structure of $\g_{0,h}$ can be read off
from the Satake diagram. In our case $\g_{0,h}$ is a
torus and $\dim\g_{0,h}=\rk\g-\dim\gt c<(\rk\g)/2$.
Hence $\dim\g_{0,h}<\dim\gt c$, and $\dim\g_0<\dim\g_1$ 
in view of Lemma~\ref{z2}.
\end{proof}

Another simple observation is that 
$\g_1{\times}\g_1\cong T^*(\g_1)$ and
the commutator morphism
$\psi: \g_1 \umt\g_1 \to \g_0$ is the moment map 
for the natural Hamiltonian action of $G_0$
on $T^*(\g_1)$. 
Therefore $\psi$ is dominant
if and only if the
generic stabiliser for the diagonal action of $G_0$ on $\gt
g_1\umt\g_1$ is finite. This is always true if the
centraliser of $\gt c$ in $\g_0$ is a torus, i.e., the
Satake diagram contains no black nodes.

\begin{lm}             \label{dim2} 
Suppose $\dim\g_1>\dim\g_0$, $\psi$ is dominant,
and $\gt C(\g_{1,h})$ is irreducible 
for all non-zero semisimple $h\in\gt
g_1$. Then $\gt C(\g_1)$ is irreducible.
\end{lm}
\begin{proof} The dominance of $\psi$ implies that if
$\psi^{-1}(\xi)\ne\varnothing$, then all irreducible components of
$\psi^{-1}(\xi)$ are of dimension $\ge 2\dim\g_1-\dim\g_0$.
Since $2\dim\g_1-\dim\g_0>\dim\g_1$, all irreducible components
of $\psi^{-1}(0)=\cgo$ are of dimension $> \dim\g_1$. 

Let $x\in\g_1$. Then $\overline{G_0(\{x\},\g_{1,x})}$
is an irreducible subvariety of $\gt C(\g_1)$ of
dimension $\dim\g_1$. If $x$ is not nilpotent, then
$\overline{G_0(\{x\},\g_{1,x})}\subset\gt C_0$, see proof
of Theorem~\ref{main}, part (1). As is well known, $G_0$ has
finitely many nilpotent orbits on $\g_1$. Hence $\gt C(\g_1)$ 
is a union of $\gt C_0$ and finite many
irreducible varieties of dimension $\dim\g_1$. Since each
irreducible component of $\gt C(\g_1)$ has dimension
greater than $\dim\g_1$,  $\gt C(\g_1)$ must be irreducible.
\end{proof}

 
\begin{thm}\label{dim3} The commuting varieties associated
with $(\gt{sl}_{2n},\gt{sl}_n{\oplus}\gt{sl}_n{\oplus}{\gt t}_1)$
and $({\EuScript E}_6,\gt{sl}_6\oplus\gt{sl}_2)$ are irreducible.
\end{thm}
\begin{proof}
First consider the classical case. 
We argue by induction on $n$. If $n=1$, then $(\g,\gt
g_0)=(\gt{so}_3,\gt{so}_2)$ and the corresponding commuting
variety is irreducible. In order to make an induction step, 
we use Lemma~\ref{dim2}. Here $\dim\g_1=2n^2$,
$\dim\g_0=2n^2-1$. The generic stabiliser for the action
of $G_0$ on $\g_1$ is an $(n{-}1)$-dimensional torus. 
It remains to show that $\gt C(\g_{1,h})$ is irreducible for all
non-zero semisimple $h\in\g_1$.

Let $h\in\g_1$ be semisimple. Using either
Lemma~\ref{sat} or the same argument as in Lemma~\ref{ss},
one can show that the semisimple part of $(\g_h,\g_{0,h})$ is a direct sum of
$\bigoplus_{i=1}^{r}(\gt{sl}_{k_i}{\oplus}\gt{sl}_{k_i},
\gt{sl}_{k_i})$ and
$(\gt{sl}_{2n-2k},\gt{sl}_{n-k}{\oplus}\gt{sl}_{n-k}{\oplus} {\gt t}_1)$.
Therefore $\gt C(\g_{1,h})$ is irreducible and we are done.

Consider now the exceptional  case. The Satake diagram of
$({\EuScript E}_6,\gt{sl}_6\oplus\gt{sl}_2)$ is:

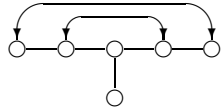
\begin{center}
\begin{picture}(70,30)(0,10)
\setlength{\unitlength}{0.017in} \put(35,3){\circle{5}}
\multiput(5,18)(15,0){5}{\circle{5}}
\multiput(8,18)(15,0){4}{\line(1,0){9}}\put(35,6){\line(0,1){9}}
\put(35,23){\oval(60,18)[t]}\put(35,23){\oval(30,10)[t]}
\multiput(5,23)(15,0){2}{\vector(0,-1){2}}
\multiput(50,23)(15,0){2}{\vector(0,-1){2}}
\end{picture}
\end{center}
It contains no black nodes,
hence $\psi$ is dominant. There are two nodes which are
not connected by an arrow with others, hence, by Lemma~\ref{dim},
$\dim\g_1>\dim\g_0$. In fact, $\dim\g_1=40$ and $\dim\g_0=38$.
To apply Lemma~\ref{dim2}, one has 
to find all Satake subdiagrams
and verify that the associated commuting varieties are
irreducible.

If  
a Satake diagram is disconnected as a graph,
then the corresponding symmetric pair is a direct sum 
of two non-trivial symmetric pairs. 
Since the commuting variety associated with a direct sum of symmetric
pairs is a product of the commuting varieties
associated with the summands, we can restrict ourselves to
connected subdiagrams. There are seven of them.

\vskip1ex
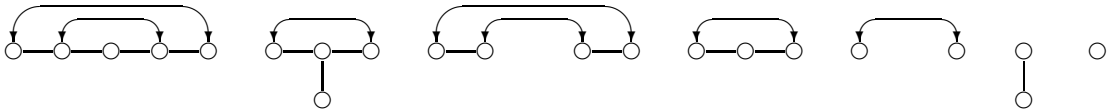
\begin{center}
\begin{tabular}{ccccccc}
\begin{picture}(70,30)(0,5)
\setlength{\unitlength}{0.017in}
\multiput(5,18)(15,0){5}{\circle{5}}
\multiput(8,18)(15,0){4}{\line(1,0){9}}
\put(35,23){\oval(60,18)[t]}\put(35,23){\oval(30,10)[t]}
\multiput(5,23)(15,0){2}{\vector(0,-1){2}}
\multiput(50,23)(15,0){2}{\vector(0,-1){2}}
\end{picture}
&
\begin{picture}(70,30)(0,5)
\setlength{\unitlength}{0.017in} \put(35,3){\circle{5}}
\multiput(20,18)(15,0){3}{\circle{5}}
\multiput(23,18)(15,0){2}{\line(1,0){9}}\put(35,6){\line(0,1){9}}
\put(35,23){\oval(30,10)[t]}
\multiput(20,23)(30,0){2}{\vector(0,-1){2}}
\end{picture}
&
\begin{picture}(70,30)(0,5)
\setlength{\unitlength}{0.017in}
\multiput(5,18)(15,0){2}{\circle{5}}
\multiput(50,18)(15,0){2}{\circle{5}}
\multiput(8,18)(45,0){2}{\line(1,0){9}}
\put(35,23){\oval(60,18)[t]}\put(35,23){\oval(30,10)[t]}
\multiput(5,23)(15,0){2}{\vector(0,-1){2}}
\multiput(50,23)(15,0){2}{\vector(0,-1){2}}
\end{picture}
&
\begin{picture}(70,30)(0,5)
\setlength{\unitlength}{0.017in}
\multiput(20,18)(15,0){3}{\circle{5}}
\multiput(23,18)(15,0){2}{\line(1,0){9}}
\put(35,23){\oval(30,10)[t]}
\multiput(20,23)(30,0){2}{\vector(0,-1){2}}
\end{picture}
&
\begin{picture}(40,30)(0,5)
\setlength{\unitlength}{0.017in}
\multiput(5,18)(30,0){2}{\circle{5}}
\put(20,23){\oval(30,10)[t]}
\multiput(5,23)(30,0){2}{\vector(0,-1){2}}
\end{picture}
&
\begin{picture}(30,30)(0,5)
\setlength{\unitlength}{0.017in} \put(15,3){\circle{5}}
\put(15,18){\circle{5}} \put(15,6){\line(0,1){9}}
\end{picture}
&
\begin{picture}(10,30)(0,5)
\setlength{\unitlength}{0.017in} \put(5,18){\circle{5}}
\end{picture}
\\
\end{tabular}
\end{center}
These diagrams correspond to the following symmetric pairs:\\
$\phantom{,}$ \quad $(\gt{sl}_6, \gt{sl}_3{\oplus}\gt{sl}_3{\oplus}{\gt t}_1)$,
$(\gt{so}_8,\gt{so}_5{\oplus}\gt{so}_3)$,
$(\gt{sl}_3{\oplus}\gt{sl}_3,\gt{sl}_3)$,\\
$\phantom{,}$ \quad  $(\gt{sl}_4, \gt{sl}_2{\oplus}\gt{sl}_2{\oplus}{\gt t}_1)=
(\gt{so}_6,\gt{so}_4{\oplus}\gt{so}_2)$,
$(\gt{sl}_2{\oplus}\gt{sl}_2,\gt{sl}_2)$,
$(\gt{sl}_3,\gt{so}_3)$, $(\gt{sl}_2,\gt{so}_2)$.

\noindent For all of them the irreducibility of $\gt C(\g_1)$ 
is already established.
\end{proof}

\section{Short $\BB Z$-gradings and reducible commuting varieties}
\label{nex}

Let $\g=\g(-1)\oplus\g(0)\oplus\g(1)$ be a
short ${\mathbb Z}$-grading of a reductive Lie algebra $\g$. Set $\g_0:=\gt
g(0)$. Then $(\g,\g_0)$ is a symmetric pair with
$\g_1=\g(-1){\oplus}\g(1)$.
Let $\gt c\subset\g_1$ be a Cartan subspace and let $\gt c(\pm 1)$ 
denote the image of $\gt c$ under the projection to $\g(\pm 1)$.
As $\gt c$ consists of semisimple elements,
we have $\gt c\cap \g(\pm 1)=0$ and $\dim\gt c=\dim\gt c(1)=\dim\gt c(-1)$.
Below, we consider the diagonal $G_0$-action on $\g(1)\times\g(1)$.

\begin{lm}\label{piki}
Suppose $\gt C(\g_1)$ is irreducible. Then
 \vskip0.7ex
\hbox to \textwidth{$(\spadesuit_1)$\hfil $\overline{G_0(\gt c(1)\umt 
\gt c(1))}=\g(1)\umt  \g(1)$. \hfil}
\end{lm}
\begin{proof}
Notice that $\g(1)$ is a commutative subalgebra, hence $\g(1)\times\g(1)\subset
\gt C(\g_1)$.
Since $\gt c\subset \gt c({-}1)\oplus \gt c(1)$, we have
\[
\overline{G_0(\gt c\umt\gt c)}\subset \overline{G_0(\gt
c(-1)\umt  \gt c(1)\umt  \gt c(-1)\umt  \gt c(1))} \subset
\overline{G_0(\gt c(-1)\umt  \gt c(-1))}\times
\overline{G_0(\gt c(1)\umt  \gt c(1))}.
\]
Hence 
$\overline{G_0(\gt c\umt\gt c)}\cap(\g(1)\umt\g(1))\subset 
\overline{G_0(\gt c(1)\umt\gt c(1))}$. If $\gt C(\g_1)$ is irreducible, then
\\
$\overline{G_0(\gt c\umt\gt c)}\cap(\g(1)\umt\g(1))=
\g(1)\times\g(1)$, and we get equality $(\spadesuit_1)$.
\end{proof}
\noindent
Similarly, one can consider condition $(\spadesuit_{-1})$ with $\gt c(-1)$
in place of $\gt c(1)$, etc. It is easily seen that
$(\spadesuit_1)$ is satisfied if and only if $(\spadesuit_{-1})$ is. 
For this reason, this common condition is denoted by $(\spadesuit)$.
\begin{cl} \label{tri} 
If condition $({\spadesuit})$ is not
satisfied, then $\gt C(\g_1)$ has at least three
irreducible components.
\end{cl}\begin{proof}
Along with the standard component $\mathfrak C_0$, we have at least
two other components
determined by $\g(1)\umt\g(1)$ and
$\g(-1)\umt\g(-1)$.
\end{proof}
\noindent
Now we give three examples of symmetric pairs arising from
short $\mathbb Z$-gradings such that condition $({\bf \spadesuit})$ is
not satisfied.

\begin{ex} \label{gl} $(\g,\g_0)=(\gt{sl}_{n+m},\gt{sl}_n\oplus\gt{sl}_m
\oplus {\gt t}_1)$.
\\
It will be easier to work with the symmetric pair
$(\gt{gl}_{n+m},\gt{gl}_n\oplus\gt{gl}_m)$ having the same commuting variety.
We use the following matrix model.
Let $V$ be an $(n+m)$-dimensional vector space and
$\g=\gt{gl}(V)$. Let $V=V_a\oplus V_b$ be a vector space decomposition
with $\dim V_a=n$ and $\dim V_b=m$. Then $\g(1)= \Hom(V_a,V_b)$,
$\g(-1)=\Hom(V_b,V_a)$, and $\g(0)=\Hom(V_a,V_a)\oplus \Hom(V_b,V_b)\simeq
\gt{gl}_n\oplus\gt{gl}_m$.
The corresponding involution $\sigma$ is the inner automorphism
determined by the matrix $A\in\GL(V)$ such that $A|_{V_a}=-\id$,
$A|_{V_b}=\id$. In the matrix form, every $\xi\in\g_1$ is depicted as follows:
$\xi=\left(\begin{array}{cc} 0 & \xi(-1) \\ \xi(1) & 0 \end{array}\right)$,
where $\xi(1)\in \Hom(V_a,V_b)$ is an $m\times n$ matrix and
$\xi(-1)\in \Hom(V_b,V_a)$ is an $n\times m$ matrix. 

Assume that $n<m$. To each pair $(\xi,\eta)\in \g_1\times\g_1$, 
we assign the $m\times 2n$ matrix
$D(\xi,\eta):=(\xi(1)\,|\, \eta(1))$. 
We will show that $\rk D(\xi,\eta)\le n$ for any $(\xi,\eta)\in 
\overline{G_0(\gt c\times\gt c)}$, whereas there are some $(\xi,\eta)\in
\g(1)\times\g(1)$ such that $\rk D(\xi,\eta)> n$. This means that 
condition $(\spadesuit)$ cannot be satisfied here.

Let $\gt c(-1)$ be the set of $n\times m$ matrices
\[
\gt c(-1)=\left\{ \left(\begin{array}{ccccccc}
        x_1 &  0     &   \ldots &  0     & 0      & \dots & 0 \\
         0  & x_2    &   \ldots &  0     & \vdots & \dots & \vdots\\
    \vdots  & \vdots &   \ddots & \vdots & \vdots & \dots & \vdots\\
        0   & \ldots &     0    & x_n    & 0      & \dots & 0 \\
        \end{array}
        \right) \mid x_1,\dots,x_n\in\bbk \right\}.
\]
Then $\gt c:=\left\{ \left(\begin{array}{cc} 
0   & M \\ 
M^t & 0 
\end{array}\right)\mid M\in \gt c(-1)\right\}$ is a Cartan subspace of $\g_1$.
It is clear from  the definition that $\rk D(\xi,\eta)\le n$ for any
$(\xi,\eta)\in\gt c\umt \gt c$. 
Next, an easy verification shows 
that this remains true for the elements of $G_0(\gt c\times\gt c)$.
Indeed, 
let $g=\begin{pmatrix} B & 0 \\ 0 & C \end{pmatrix}\in G_0$. 
Then $\Ad(g){\cdot}\xi=\begin{pmatrix} 0 & B\xi(-1)C^{-1} \\
C\xi(1)B^{-1} & 0 \end{pmatrix}$. It follows that
$D(\Ad(g)\um \xi,\Ad(g)\um\eta)= CD(\xi,\eta)\hat B$, where
$
\hat B=\begin{pmatrix}
B^{-1} & 0\\ 0 & B^{-1} \end{pmatrix}
$.
Hence $\rk D(\xi,\eta)=
\rk D(\Ad(g)\um \xi,\Ad(g)\um\eta)$ for any $g\in G_0$. 

On the other hand, it is clear that there are two $m\times n$ matrices
$\xi(1),\eta(1)$, i.e., two elements of $\g(1)$, 
such that $\rk(\xi(1)\,|\,\eta(1))=\min(2n,m)>n$, since $m>n$.
This completes our argument.

\end{ex}
\noindent
Thus, $\gt C(\g_1)$ 
has at least three irreducible components, by virtue of 
Corollary~\ref{tri}. However, this lower bound is far from
being precise. Surprisingly, $\gt C(\g_1)$ can have
arbitrary many irreducible components.

\begin{prop}         \label{irrcomp} 
Suppose $(\g,\g_0)=(\gt{gl}_{n+m},\gt{gl}_n{\oplus}\gt{gl}_m)$ with $n\le m$.
Then $\gt C(\g_1)$ has at least $2\min(2n,m)-2n+1$ irreducible components.
\end{prop}
\begin{proof} 
For $\xi=\begin{pmatrix} 0 & Y \\ X & 0 \end{pmatrix}\in\g_1$ and 
$\eta=\begin{pmatrix} 0 & U \\ Z & 0 \end{pmatrix}\in\g_1$, 
we set $D_1(\xi,\eta):=(X\,|\,Z)$
(as above) and $D_2(\xi,\eta):=(Y^t\,|\,U^t)$.
These are two $m\times 2n$-matrices.
\\[.7ex]
{\sf 1.} 
First, we are going to prove that 
$\rk D_1(\xi,\eta)+\rk D_2(\xi,\eta)\le 2n$ if\/ $[\xi,\eta]=0$.
Suppose $\rk X=r$. Replacing 
$(\xi,\eta)$ with $(\Ad(g){\cdot}\xi,\Ad(g){\cdot}\eta)$ 
for a suitable $g\in G_0$, we may assume that 
$X=(x_{ij})$ such that $x_{ii}=1$ if
$i=1,\ldots,r$ and $x_{ij}=0$ for all other $i,j$. Then 
$\rk D_1(\xi,\eta)=r+\rk Z_1$, 
where $Z_1$ is the $(m-r){\times}n$ matrix 
consisting of the last $(m-r)$ rows of $Z$. 
Condition $[\xi,\eta]=0$ means that $XU=ZY$ and 
$YZ=UX$. 
The last $r$ rows of $X$ and hence of $XU$
are zero. Therefore $Z_1Y=0$ and $\rk Y\le n-\rk Z_1$.

Let $X_0$, $Z_0$, $U_0$ be the sub-matrices of $X$, $Z$, and $U$, respectively,
consisting of the first $r$ rows. From the equality 
$XU=ZY$ and our assumption on $X$, we then obtain $Z_0Y=X_0U=U_0$. 
Therefore the rank of the $m{\times}(n+r)$ matrix $(Y^t\,|\,U_0^t)$ 
is not grater than $\rk Y$ and hence $\rk(Y^t\,|\,U^t)\le\rk Y+(n-r)$. 
Thus, 
\[
\rk D_1(\xi,\eta)+\rk D_2(\xi,\eta)\le
(r+\rk Z_1)+(n-\rk Z_1+n-r)=2n.
\]
{\sf 2}. Let $P_q$ be the closed subset of $\gt C(\gt g_1)$
defined by the conditions:
\[
\rk D_1(\xi,\eta)\le q\ \text{ and }\ \rk D_2(\xi,\eta)\le 2n-q.
\] 
As was just proved, $\gt C(\gt g_1)=\bigcup P_q$. 
Since $D_1(\xi,\eta)$ is an $m{\times}2n$-matrix, 
$\rk D_1(\xi,\eta)\le\min(2n,m)$. 
Hence $P_q\subset P_{\min(2n,m)}$ if $q>\min(2n,m)$. 
Similarly, for $q<2n-\min(2n,m)$, we have $P_q\subset P_r$, where 
$r={2n-\min(2n,m)}$. This implies that
\[
\gt C(\gt g_1)=\displaystyle \bigcup_{2n-\min(2n,m)}^{\min(2n,m)} P_q \ .
\]
Let us show that $P_q\not\subset\bigcup\limits_{r\ne q} P_r$ 
for each $q\in [2n-\min(2n,m),\min(2n,m)]$. Since 
$\gt C(\g_1)$ is invariant under the transpose, 
there is no harm in assuming that $q\ge n$.

Keep the previous notation. Take 
$(\xi,\eta)\in\gt C(\g_1)$ such that the submatrices $X_0, Z_0, U_0$, etc. 
satisfy the following conditions: 
$\rk X=\rk X_0=n$ (i.e., $r=n$), $U=0$, $\rk Z_1=q-n$, and $Z_0=0$. 
Let $Y_0$ and $Y_1$ denote the 
submatrices of $Y$ consisting of the first $n$ columns and 
the last $m-n$ columns, respectively. Then the
condition $[\xi,\eta]=0$ is equivalent to that $Z_1Y=0$ and $Y_1Z_1=0$. 
One can choose $Y$ satisfying 
these conditions and such that $Y_1=0$ and $\rk Y_0=n-\rk Z_1$. 
Then $\rk D_1(\xi,\eta)=q$ and $\rk D_2(\xi,\eta)=2n-q$.
That is, we found a point in 
$P_q\setminus\bigcup\limits_{r\ne q} P_r$.
It follows that each $P_q$ contains an irreducible component of
$\gt C(\g_1)$ that does not belong to the other $P_r$, and we are done.
\end{proof}
\noindent
{\bf Remark.} Motivated by the formula of Proposition~\ref{irrcomp}, we set
$F(n,m)=2\min(2n,m)-2n+1$ if $1\le n\le m$.
For $n=1$ and $m\ge 2$, we have 
$\rk(\gt{gl}_{m+1},\gt{gl}_m\oplus\gt{gl}_1)=1$, 
and it was shown in \cite{Dima2} that $\cgo$ has three irreducible
components. Since $F(1,m)=3$ if $m\ge 2$, Proposition~\ref{irrcomp} gives here 
the exact value for the number of irreducible components.
For $m=n$, we also obtain the exact number of irreducible components,
in view of Theorem~\ref{dim3}.
It is curious to observe that $F(n,m)=F(m-n,m)$.
Hence $F(m-1,m)=3$, and one may ask whether it is true
that, for the symmetric pair 
$(\gt{gl}_{2m-1},\gt{gl}_m\oplus\gt{gl}_{m-1})$, the variety
$\cgo$ has exactly three irreducible components. More generally, one
may conjecture that the varieties $\{P_q\}$ in the previous proof are 
irreducible and therefore the number of irreducible components is 
always equal to $F(n,m)$.
The equality $F(n,m)=F(m-n,m)$ also suggests that there might be a natural
bijection between the irreducible components of the commuting varieties
associated with $(\gt{gl}_{m+n},\gt{gl}_m\oplus\gt{gl}_{n})$ and
$(\gt{gl}_{2m-n},\gt{gl}_m\oplus\gt{gl}_{m-n})$.

\begin{ex} \label{so} $(\g,\g_0)=(\gt{so}_{2n},\gt{gl}_n)$.
\\
Let $V$ be a $2n$-dimensional vector space and 
$\g=\gt{so}(V)=\gt{so}_{2n}$.
Consider a decomposition of $V$ into a
direct sum of two isotropic subspaces $V=V_+\oplus V_-$.
Let $A\in{\rm O}(V)$ be such that $A|_{V_+}=\id$ and
$A|_{V_-}=-\id$. Then the conjugation by $A$ defines an
involution of $\g$ such that 
$\g_0=\gt{gl}_n\cong\gt{gl}(V_+)\cong\gt{gl}(V_-)$ and
$\g_1\subset\Hom(V_+,V_-)\oplus\Hom(V_-,V_+)$. More
precisely, one can choose  bases for $V_+$ and $V_-$ such
that
\[
G_0=\left\{\begin{pmatrix}
          B & 0 \\
          0 & (B^t)^{-1} 
         \end{pmatrix}\mid B\in\GL_n \right\}, \enskip
\g_1=\left\{\begin{pmatrix}
         0 & X \\
         Y & 0 \\
         \end{pmatrix} \mid X,Y\in \gt{gl}_n,
 X=-X^t, Y=-Y^t \right\}.
\]
Here $\g(1)=\g_1\cap\Hom(V_+,V_-)$ and
$\g(-1)=\g_1\cap\Hom(V_-,V_+)$.
For a pair $(\xi,\eta)\in\g_1\times\g_1$ with
$\xi=\begin{pmatrix}
         0 & X \\
         Y & 0 \\
         \end{pmatrix}$ and $\eta=
    \begin{pmatrix}
         0 & Z \\
         U & 0 \\
         \end{pmatrix}$,
we set $D_1(\xi,\eta):=(X\,|\, Z)$. 
We take the Cartan subspace $\gt c\subset\g_1$ consisting of
skew-symmetric anti-diagonal matrices, i.e.,
\[
\gt c=\left\{
\begin{pmatrix}
         0 & X \\
         X & 0 \\
\end{pmatrix} \mid  X=-X^t \textrm{ and } X=(x_{ij}),\ x_{ij}=0 \textrm{ if }
i+j\ne n+1 \right\}.
\]
Now assume that $n=2k+1$. Then $\rk D_1(t,h)<n$
for any $(t,h)\in\gt c\times\gt c$. If
$g=\begin{pmatrix}
          B & 0 \\
          0 & (B^t)^{-1} \\
         \end{pmatrix}\in G_0$
and $\xi=\begin{pmatrix}
         0 & X \\
         Y & 0 \\
         \end{pmatrix}\in\g_1$, then
$\Ad(g)\um \xi=\begin{pmatrix}
         0 & BXB^t \\
         (B^t)^{-1}YB^{-1} & 0 \\
         \end{pmatrix}$. Therefore,
$
D_1(\Ad(g)\um \xi,\Ad(g)\um \eta)=
BD_1(\xi,\eta)\hat B$, where 
$\hat B=\begin{pmatrix}
        B^t & 0     \\
          0 & B^{t} \\
\end{pmatrix}$ 
is a non-degenerate $2n\umt 2n$ matrix.  Hence
$\rk D_1(\Ad(g)\um \xi,\Ad(g)\um \eta)=\rk D_1(\xi,\eta)$
and $\rk D_1(t,h)< n$ for each pair $(t,h)\in\gt C_0$.

Let $X,Z\in\g(1)$ be skew-symmetric $n\umt n$ matrices
of rank $2k$ such that the last column and the last row of
$X$ are zero, and the first row and the first column of $Z$
are zero. Clearly, if $k\ge 1$, then $\rk(X|Z)=2k+1=n$ 
and therefore $(X,Y)\not\in\overline{G_0(\gt c(1)\umt  \gt
c(1))}$. Thus, for $(\gt{so}_{2n},\gt{gl}_n)$ with odd
$n\ge 3$, condition $({\bf \spadesuit})$ is not satisfied.
\end{ex}

\begin{ex} \label{e6k}
Consider now the symmetric pair
$({\EuScript E}_6,\gt{so}_{10}\oplus {\gt t}_1)$. 
\\
Here $\g(1)$ and $\g(-1)$ are different half-spinor 
representations of $\gt{so}_{10}$. Let $\gt t$ be a Cartan 
subalgebra of $\gt{so}_{10}$ and $\{\esi_1,\ldots, \esi_5\}$ 
an orthonormal basis of $\gt t^*_{\BB Q}$.
Let  $\{\pi_1,\ldots, \pi_5\}$  be
the fundamental weights of $\gt{so}_{10}$
such that
\[
\pi_4=(\esi_1+\esi_2+\esi_3+\esi_4-\esi_5)/2 \ 
\textrm{ and }\ 
\pi_5=(\esi_1+\esi_2+\esi_3+\esi_4+\esi_5)/2 \ 
\]
are the highest weights of half-spinor representations.
(See \cite[Reference Chapter]{vo} for more details on the notation.) 
Let $\mathcal R(\varphi)$ denote the simple $\gt{so}_{10}$-module with
highest weight $\vp$. We assume that $\g(1)=\mathcal R(\pi_5)$ and
$\g(-1)=\mathcal R(\pi_4)$. The rank of
$({\EuScript E}_6,\gt{so}_{10}\oplus {\gt t}_1)$ equals $2$ and a Cartan
subspace $\gt c\subset\g_1$ can be chosen such that $\gt
c(1)$ is a $\gt t$-stable subspace of $\g(1)$ with weights $\pi_5$
and $(\esi_1-\esi_2-\esi_3-\esi_4-\esi_5)/2$. 
\\[.6ex]
Because the weights of $\gt c(1)$ belong to an open halfspace of 
$\gt t^*_{\BB Q}$, $\ce(1)\times \ce(1)$ lies in the nullcone of the
$\Spin_{10}$-module $\g(1)\times\g(1)$. (Actually, similar assertion 
holds for any number of copies of $\ce(1)$ and $\g(1)$.
Therefore if $f\in \bbk[\g(1)\umt\g(1)]^{\Spin_{10}}$ is a homogeneous
polynomial of positive degree, then it vanishes on 
$\overline{G_0{\cdot}(\gt c(1)\umt  \gt c(1))}$. Hence the very existence of 
non-trivial $\Spin_{10}$-invariant polynomials on $\g(1)\umt\g(1)$
will imply that $(\spadesuit)$ is not satisfied.

As ${\mathcal S}^2(\mathcal R(\pi_4))=\mathcal R(2\pi_4)\oplus 
\mathcal R(\pi_1)$, there are
the following  $\Spin_{10}$-equivariant inclusions:
\[
{\mathcal S^4}(\mathcal R(\pi_4){\oplus}\mathcal R(\pi_4))\supset
{\mathcal S^2}(\mathcal R(\pi_4))\otimes{\mathcal S^2}(\mathcal R(\pi_4))\supset
 \mathcal R(\pi_1)\otimes \mathcal R(\pi_1).
\]
Since $(\mathcal R(\pi_1){\otimes}\mathcal R(\pi_1))^{\Spin_{10}}\ne 0$, 
we obtain a desired
non-trivial ${\Spin}_{10}$-invariant of degree 4 in
$\bbk[\g(1){\times}\g(1)]=
{\mathcal S^\bullet}(\mathcal R(\pi_4){\oplus}\mathcal R(\pi_4))$.
\end{ex}

\section{Principal nilpotent pairs and commuting varieties}
\label{pnp}

\noindent
In this section, we prove that $\cgo$ is reducible for two symmetric pairs
related to exceptional Lie algebras.

Let $\gt s$ be the centraliser in $\g$ of a generic element of $\ce$.
Since $\ce$ is diagonalisable, 
$\g_x\cap\g_y \supset\gt s$ for any $x,y\in\ce$.
Therefore, $\dim(\g_x\cap\g_y) \ge \dim\gt s$ for any pair
$(x,y)\in \overline{G_0(\ce\times\ce)}=\gt C_0$. 
In what follows, we write $\g_{x,y}$ for $\g_x\cap\g_y$.
Assume that we managed to find a pair of commuting elements 
$\xi,\eta$ in $\g_1$ such that 
$\dim\g_{\xi,\eta}< \dim\gt s$. Then clearly
$(\xi,\eta)\not\in\gt C_0$, and we conclude that $\cgo$ is reducible.

Recall that $\dim\g_{\xi,\eta}\ge\rk\g$ for any pair of commuting elements.
This follows from the irreducibility of the "usual" commuting variety
$\gt C(\g)$. Therefore, one may try to find a symmetric pair such that
$\dim\es > \rk\g$, whereas $\dim\g_{\xi,\eta}=\rk\g$
for some $(\xi,\eta)\in\cgo$. This will be done with the help of theory
of principal nilpotent pairs. 

\begin{df}[Ginzburg~\cite{vit}]
A pair $(e_1,e_2) \in \g\times\g$
is said to be {\it principal nilpotent\/} 
(\pn pair, for short) if

\centerline{$[e_1, e_2]=0$, \quad $\dim\g_{e_1,e_2}=\rk\g$,} 

\noindent 
and there exists
a pair of semisimple elements $(h_1, h_2)\in\g\times\g$ such that 
$\ad h_1$ and $\ad h_2$ have integral eigenvalues and
$[h_1, h_2]=0,\quad [h_i, e_j]=\delta_{ij} e_j \quad
(i,j\in \{1,2\})$.  
\end{df}

Since the eigenvalues of $\ad h_1$ and $\ad h_2$ are integral,
we can consider the corresponding $\BB Z\times\BB Z$-grading of $\g$:
\[
  \g=\bigoplus\g(i,j), \text{ where } \g(i,j)=
\{x\in\g\mid [h_1,x]=ix, [h_2,x]=jx \}.
\]
Notice that $e_1\in\g(1,0)$ and $e_2\in\g(0,1)$. Using this
bi-grading, we define a symmetric pair ($\BB Z_2$-grading) by letting
$\displaystyle \g_0=\bigoplus_{i+j \text{ is even}} \g(i,j)$.
From this construction, it follows that $(e_1,e_2)\in\cgo$.
If we are lucky and the resulting symmetric pair has the property that
$\dim\es > \rk \g$, i.e., $\es$ is not Abelian, 
then we certainly know that $\cgo$ is reducible.

\noindent 
It was shown by Ginzburg \cite{vit} that the number of $G$-orbits of 
\pn pairs is finite,
and a complete classification of such pairs in exceptional Lie 
algebras is obtained by Elashvili and Panyushev \cite[Appendix]{vit}.
Using that classification, it is not hard to find out which \pn pairs 
lead to symmetric pairs (as described above) such that
$\dim\es > \rk \g$. 
Actually, there is a unique \pn pair with this property if $\g$
is exceptional. This pair occurs in the next theorem.

\begin{thm}  \label{e7+e8}  
For the symmetric pairs $({\EuScript E}_7,\gt{so}_{12}{\oplus}\gt{sl}_2)$ 
and $({\EuScript E}_8,{\EuScript E}_7{\oplus}\gt{sl}_2)$, 
the variety $\cgo$ is reducible.
\end{thm}\begin{proof}
1. For the simple Lie algebra ${\EuScript E}_7$, there is a \pn 
pair $(e_1,e_2)$ such that both $e_1,e_2$ belong to the $E_7$-orbit
of type $\eus A_4+\eus A_1$. (See below an explicit construction of such a
pair.)
The bi-grading corresponding to this pair is depicted
in \cite[Fig.\,2]{genpp}. This shows that $\dim\g_0=69$ for the corresponding
$\BB Z_2$-grading. Hence $\g_0\simeq \gt{so}_{12}{\oplus}\gt{sl}_2$.
For this symmetric pair, we have $[\es,\es]\simeq 3{\eus A}_1$. Hence
$(e_1,e_2)\not\in\gt C_0$. 
\\[.7ex]
2. The same \pn pair can be used for proving reducibility 
in the $\EuScript E_8$-case.
Using the standard Lie algebra embedding 
$\g:={\EuScript E}_7 \subset {\EuScript E}_8=:\tilde\g$,
the first symmetric pair can be regarded as a subpair of the second.
This simply means that we have two embeddings
$\g_0=\gt{so}_{12}{\oplus}\gt{sl}_2\subset 
{\EuScript E}_7{\oplus}\gt{sl}_2=\tilde\g_0$
and $\g_1\subset\tilde\g_1$. Then $(e_1,e_2)\in\cgo\subset \gt C(\tilde\g_1)$.
As above, we wish to prove that $(e_1,e_2)$ does not  belong to the
irreducible component 
$\tilde{\gt C}_0=\overline{\tilde G_0(\tilde\ce\times\tilde\ce)}$.
The two symmetric pairs in question have the
same rank (namely, 4), so we actually may assume that $\tilde\ce=\ce$.
For a generic element of $\ce$, its stabiliser in $\tilde\g$ has the 
semisimple part of type $\eus D_4$. Therefore $\dim\tilde\es=28+4=32$. 
On the other hand, we will prove that 
$\dim\tilde\g_{e_1,e_2}=26< 32$, which establishes the 
desired reducibility.

To this end, we need an explicit construction of the \pn pair
$(e_1,e_2)$ in ${\EuScript E}_7$. 
According to classical results of E.B.\,Dynkin~\cite{dy}, 
${\EuScript E}_7$ contains a unique
(up to conjugation) maximal semisimple subalgebra of type $\eus A_2$.
The embedding $\eus A_2\subset {\eus E}_7$ has many remarkable properties.
For instance, ${E}_7/A_2$ is an isotropy irreducible homogeneous space.
Here ${\eus E}_7/\eus A_2$ is the simple $\eus A_2$-module 
$\mathcal R(4\pi_1{+}4\pi_2)$. 
Let $\ap,\beta,\ap+\beta$ be the positive roots of $\eus A_2$.
There are two non-equivalent choices of nontrivial \pn pairs in
$\eus A_2$. One takes $(e_1,e_2)=(e_\ap,e_{\ap+\beta})$ or 
$(e_\beta,e_{\ap+\beta})$. 
Then one easily verifies that these two pairs remain principal also in
$\eus E_7$. Furthermore, one can prove that these two pairs are conjugate
with respect to the group $E_7$, so one obtains a single conjugacy class
of \pn pairs in $\eus E_7$.

In order to compute $\dim\tilde \g_{e_1,e_2}$, we first decompose 
$\tilde\g=\eus E_8$ as $\eus E_7$-module, and then further as $\eus A_2$-module.
First we notice that
\[
   \tilde\g\vert_{\eus E_7}=\g \oplus 2\mathcal V \oplus 3\odin \ ,
\]
where $\mathcal V$ is the 56-dimensional simple $\eus E_7$-module
and $\odin$ is the trivial one-dimensional module. Since $(e_1,e_2)$
is a \pn pair in $\eus E_7$, we have $\dim \g_{e_1,e_2}=7$.
To compute the fixed-point subspace of $(e_1,e_2)$ in $\mathcal V$,
we use fact that $\mathcal V\vert_{\eus A_2}=\mathcal R(6\pi_1)\oplus
\mathcal R(6\pi_2)$, see \cite[Table 24]{dy}. 
Finally, one easily computes using the above two presentations of $(e_1,e_2)$
that $\dim \mathcal R(6\pi_1)^{e_1,e_2}=1$ and 
$\dim \mathcal R(6\pi_2)^{e_1,e_2}=7$ (or vice versa).
Altogether, we obtain
\[
   \dim\tilde \g_{e_1,e_2}=7+2(7+1)+3=26.
\]
This completes the proof of reducibility in the $\eus E_8$-case.
\end{proof}

\section{A few remarks on the remaining symmetric pairs}
\label{unk}

\noindent
There are three cases
for which the irreducibility of $\gt C(\g_1)$ is not known yet.
They are listed below, together with their Satake diagrams.

\begin{tabular}{ccc}
1) & $(\gt{sp}_{2n+2m},\gt{sp}_{2n}{\oplus}\gt{sp}_{2m})$,
$n\ge m\ge 2$ & 
\begin{picture}(150,25)(0,5)
\multiput(30,8)(60,0){2}{\circle{6}}
\multiput(70,8)(40,0){3}{\circle*{6}}
\multiput(10,8)(160,0){2}{\circle*{6}}
\multiput(152.5,7)(0,2){2}{\line(1,0){15}}
\multiput(73,8)(20,0){2}{\line(1,0){14}}
\multiput(33.1,8)(29,0){2}{\line(1,0){5}}
\multiput(113,8)(29,0){2}{\line(1,0){5}}
\put(42,5){$\cdots$}
\put(122,5){$\cdots$}
\put(13,8){\line(1,0){14}}
\put(8,0){$\underbrace%
{\mbox{\hspace{86\unitlength}}}_{2m}$}
\put(155,5){$<$} 
\end{picture}    \\
2) & $(\gt{so}_{4n},\gt{gl}_{2n})$, $n\ge 2$ & 
\begin{picture}(150,30)(-10,5)
\multiput(10,8)(40,0){2}{\circle*{6}}
\multiput(30,8)(40,0){2}{\circle{6}}
\multiput(13,8)(20,0){3}{\line(1,0){14}}
\put(110,8){\circle*{6}} \put(130,8){\circle{6}}
\put(113,8){\line(1,0){14}}
\multiput(150,-2)(0,20){2}{\circle{6}}
\put(150,18){\circle*{6}} \put(133,10){\line(2,1){13}}
\put(133,6){\line(2,-1){13}}
\multiput(73,8)(29,0){2}{\line(1,0){5}}
\put(82,5){$\cdots$}
\end{picture}                    \\
3) & $({\EuScript E}_7, {\EuScript E}_6{\oplus}{\gt t}_1)$ &
\begin{picture}(90,30)(5,7)
\setlength{\unitlength}{0.016in} \multiput(23,18)(15,0){5}{\line(1,0){9}}
\put(65,3){\circle*{5}}\multiput(20,18)(15,0){3}{\circle{5}}\put(65,6){\line(0,1){9}}
\multiput(50,18)(15,0){3}{\circle*{5}}
\put(95,18){\circle{5}}\end{picture} \\
%
\end{tabular}
\vskip.5ex
\noindent
Items~2) and 3) are symmetric pairs of Hermitian type, i.e.,
$\g_0$ is not semisimple. Such symmetric pairs (and the corresponding
symmetric spaces) can be either of {\it tube} or {\it non-tube type}. 
All symmetric pairs of Hermitian type arise form short $\mathbb Z$-gradings
of $\g$. This construction is explained in Section~\ref{nex}.  
There are many characterisations
(definitions) of symmetric pairs of tube type. 
One of the possible definitions is the following.
Let $G_0'$ be the commutator group of $G_0$. 
Then $(\g,\g_0)$ is said to be of tube type if 
$\bbk[\g(1)]^{G_0'}\ne \bbk$. 
Then one can prove that
$(\g,\g_0)$ is of tube type if and only if the dense (nilpotent)
orbit in $G{\cdot}\g(1)$ is even  if and only if the $G_0'$-action on
$\g(1)$ is stable. Yet another characterisation is that the symmetric 
pairs of tube type are those related to simple Jordan algebras. 

Examples~\ref{gl} (with $n\ne m$), \ref{so} (with $n$ odd), and \ref{e6k}
represent all Hermitian symmetric pairs of non-tube type, and we have 
proved that for all these cases $\gt C(\g_1)$ is reducible. On the other hand,
for the symmetric pair $(\gt{sl}_{2n},\gt{sl}_n\oplus\gt{sl}_n\oplus {\gt t}_1)$, which
is of tube type, $\gt C(\g_1)$ is irreducible, as proved in Section~\ref{e6}.
Other symmetric pairs of tube type are 
$(\gt{sp}_{2n}, \gt{gl}_n)$ and $(\gt{so}_{n+2}, \gt{so}_n\oplus\gt{so}_2)$,
where the irreducibility is also known.
This suggests that for items~2) and 3), which are of tube type,
$\gt C(\g_1)$ ought to be irreducible, too.
Hopefully, the connection with Jordan algebras
might lead to a uniform proof of the irreducibility. 

The first case is probably the most intricate. 
Here $\rk(\gt{sp}_{2n+2m},\gt{sp}_{2n}{\oplus}\gt{sp}_{2m})=\min(n,m)=m$.
Since the commuting variety associated with 
$(\gt{sp}_{2n+2},\gt{sp}_{2n}{\oplus}\gt{sp}_{2})$
is reducible, there is a non-zero semisimple 
element $h\in\g_1$ such that $\gt C(\g_{1,h})$ is reducible. 
Therefore, our induction scheme does not work in this situation.
Nevertheless, using {\sl ad hoc} methods the second author has proved
that $\gt C(\g_1)$ is irreducible for
$(\gt{sp}_{2n+4},\gt{sp}_{2n}{\oplus}\gt{sp}_{4})$, $n\ge 2$.

The problem of irreducibility of $\gt C(\g_{1,h})$ also
arises for second and third cases. The second case might be treated by 
induction, since all proper sub-symmetric pairs are direct sums of
pairs of the same form and pairs $(\gt{sl}_{2m},\gt{sp}_{2m})$, where 
the irreducibility is known. For the third case,
all proper sub-symmetric pairs $(\g_h, \g_{0,h})$ can be found 
using Proposition~\ref{sat}. 
It turns out that here 
all commuting varieties $\gt C(\g_{1,h})$ are irreducible. 

Another question is whether a nilpotent 
$G_0$-orbit in $\g_1$ gives rise to an irreducible component
of $\gt C(\g_1)$ different from $\gt C_0$. 
Provided all $\gt C(\g_{1,h})$ are irreducible, 
this question can be reduced to non-even
$\sigma$-distinguished nilpotent elements. 
A classification of $\sigma$-distinguished nilpotent 
orbits is deduced from Djokovi{\'c}'s calculations \cite{dj}. 
According to Table~13 of~\cite{dj}, for item~3),
there are two non-even $\sigma$-distinguished 
$G_0$-orbits in $\g_1$. 
Thus, the task of verifying irreducibility of 
$\gt C(\g_1)$ in the third case seems to
be manageable.

\end{document}